\documentclass[a4paper,12pt,legno]{article}

\usepackage{amsmath}
\usepackage{amsfonts}
\usepackage{amsthm}
\usepackage{mathrsfs}
\usepackage{amssymb}

\newtheorem{thm}{THEOREM}
\newtheorem{lem}{LEMMA}
\newtheorem{prop}{PROPOSITION}
\newtheorem{cor}{COROLLARY}
\newtheorem{ex}{Example}

\theoremstyle{remark}
\newtheorem{rem}{\textbf{Remark}}

\usepackage{graphicx}
\usepackage[T1]{fontenc}
\usepackage[polish]{babel}
\usepackage[utf8]{inputenc}

\newcommand{\Rset}{\mathbb{R}}
\newcommand{\Cset}{\mathbb{C}}

\newcommand{\al}{\alpha}
\newcommand{\be}{\beta}

\begin{document}

\title{Calculus on random integral mappings $I^{h,r}_{(a,b]}$ and their domains\footnote{Research funded by Narodowe Centrum Nauki (NCN)
Dec2011/01/B/ST1/01257}}

\author{ Zbigniew J. Jurek (University of Wroc\l aw\footnote{Part of this work was done when Author was visiting
Indiana University, Bloomington, \qquad \qquad \qquad \qquad USA in
Spring 2013.}\,\,)}

\date{October 2 , 2013}

\maketitle
\begin{quote} \textbf{Abstract.} It is proved that the random integral
mappings (some type of functionals of L\'evy processes) are 
isomorphisms between convolution semigroups of infinitely divisible
measures. However, the inverse mappings are no longer of the random
integral form. Domains are characterized in many ways. Compositions
(iterated integrals) can be expressed as a single random integral
mapping. Finally,  obtained results are illustrated by examples.



\emph{Abbreviated title: Calculus on random integral mappings}

\end{quote}
\maketitle

Addresses:

\medskip
\noindent Institute of Mathematics \\ University of Wroc\l aw  \\
Pl. Grunwaldzki 2/4
\\ 50-384 Wroc\l aw \\ Poland. \\ www.math.uni.wroc.pl/$\sim$zjjurek ;
e-mail: zjjurek@math.uni.wroc.pl

\newpage
For an interval $(a,b]$ in the positive half-line, two deterministic
functions $h$ and $r$,  and a L\'evy process $Y_{\nu}(t), t\ge0$,
where $\nu$ is the law of random variable $Y_{\nu}(1)$, we consider
the following mapping
\[
\nu\longmapsto
I^{h,r}_{(a,b]}(\nu):=\mathcal{L}\big(\int_{(a,b]}h(t)\,dY_{\nu}(r(t))\big),
\ \ (\star)
\]
where $\mathcal{L}$ denotes the probability distribution of the random
(stochastic) integral. One of the problems related to $(\ast)$ is to
show that the mappings $I^{h,r}_{(a,b]}$  are one-to-one  and  to
characterize their domains. We consider here both questions for
fairly general classes of functions $h$ and $r$ and measures $\nu$ (
L\'evy processes $Y_{\nu}$) on a real separable Banach space.

Let us recall that over the past decades the method of describing 
a given measure as a probability distribution of an integral
$(\ast)$ was successfully applied in many instances. Already  in
Jurek-Vervaat (1983) it was proved that in order that
\[
a_n(\xi_1+\xi_2+...+\xi_n)+x_n\Rightarrow \mu,   \ \ \  (\star
\star)
\]
for some infinitesimal triangular array $a_n\xi_j, \, 1\le j \le
n; \ n\ge 1$, it is necessary and sufficient that
\[
\mu =\mathcal{L}(\int_{(0, \infty)}\,e^{-t}dY_{\nu}(t))\equiv
I^{e^{-t},\,t}_{(0,\infty)}(\nu) , \ \ \ \ \ (\star \star \star)
\]
where the L\'evy process $(Y_{\nu}(t), 0 \le t < \infty)$ is such
that $\nu$ has finite logarithmic moment.  The expression $(\star \star \star)$ was called \emph{a random
integral representation} of the selfdecomposable (or L\'evy class
$L$) measure $\mu$ and $Y_{\nu}$ was referred to as \emph{the
background driving L\'evy process} (BDLP) of $\mu$.

The phenomena of identifying a class of limit distributions  with a
collection of laws of random integrals
 $(\ast)$  was proved for many other
limiting schemes. In Jurek (1988),(1989) and Jurek-Schreiber (1992)
almost the whole class ID, of all infinitely divisible probability
measures, was described as a sum of increasing subsemigroups.
More precisely, if
\[
\mathcal{U}_{\beta}:=\{
I^{t,t^{\beta}}_{(0,1]}(\nu)=\mathcal{L}(\int_{(0,1]}t\,dY_{\nu}(t^{\beta})):\nu\in
ID\},
\]
then
\[
ID=\overline{ ( \cup_{\beta_n}\mathcal{U}_{\beta_n})}
\]
\noindent for any sequence $\beta_n$ increasing  to infinity and the bar means closure in the  weak topology.

Still later on, many new classes of  probability distributions were
simply defined as laws of some random integrals analogous to $(
\star)$  without any reference to limiting procedures. We illustrate that approach with two papers.
Sato (2006) introduced two specific families of random
integrals on $\Rset^d$ by specifying the inner clock $r$ in
$(\star)$. One of them had the time change function
\[
r(t):=\int_t^{\infty}u^{-\alpha-1}\,e^{-u}\,du,
\]
and the space scaling $h(t)=t$. On the other hand, in Maejima,
Perez-Abreu and Sato (2012), the authors introduced subclasses of
infinitely measures by specifying the map $(\ast)$ not in terms of
measures $\nu$ but in terms of their L\'evy (spectral) measures M;
comp. formula 10 (iii) below. Using the arcsine density for the time
change $r$, the authors introduced two transforms $\mathcal{A}_{1}$ and
$\mathcal{A}_{2}$ and defined the corresponding subclasses of infinitely probability
measures. The one given by $\mathcal{A}_{2}$ gives the mapping $(\ast)$ with $h(t)=t$ and $r$ as the cumulative
distribution function of the arcsine density on $(0,1)$.

In this paper we propose a quite general approach to random integral
representations and mappings. For a more convenient way of
navigating the body of research in random representations and an
ease of comparing results of different authors we propose here a new
and unified form of definitions and notations.

Finally, we will utilize here the notion and properties of image
measures, in particular, images under the tensor product of
functions. Our results are new  on  $\Rset^d$ (Euclidean space). Most of  proofs are given 
in the generality of measures on  real separable Banach spaces. However, no essential knowledge of the functional analysis
is required.

Last but not least, the proposed calculus on random integral
mappings and their domains might be formally viewed as an analogue
of  the calculus on linear operators on Banach spaces and their
domains (in functional analysis).

\medskip
\medskip
\textbf{0. Notations and brief descriptions of results.}

\medskip
\emph{0.1. \underline{Notations and basic facts}}

$E$ is a real separable Banach space;

$E^{\prime}$ is the topological dual of $E$;

$<\cdot,\cdot>$ is the dual pair (scalar product) between
$E^{\prime}$ and $E$;

$\Rightarrow$ denotes the weak convergence of probability measures;

$\mathcal{L}(X)$ is the probability distribution of random variable
$X$;

$\stackrel{d}{=}$ means equality in distribution;

$D_F[a,b]$ denotes the Skorochod space of F-valued functions that
are right

continuous on $[a,b)$ with left-hand limits on $(a,b]$ ; in short:
\emph{cadlag}

functions;  $F$ is a complete separable metric space;

$(Y_{\nu}(t), 0\le t < \infty)$  denotes a L\'evy process such that
$\mathcal{L}(Y_{\nu}(1))=\nu$;

$ID(E)$ denotes the set of all infinitely divisible Borel measures
on $E$;

$\hat{\mu}$ is the characteristic functional (Fourier transform) of
$\mu$;

$\Phi(y)=\log\hat{\mu}(y)$ is the \emph{L\'evy exponent} of $\mu \in
ID(E)$, ($y \in E^{\prime}$);

$\Phi(y)=i<y,z_0>-\frac{1}{2}<y,Ry>$

$+ \int_{E\setminus{\{0\}}}(e^{i<y,x>}-1-i<y,x>1_B(x))M(dx)$
\emph{(L\'evy-Khintchine formula)};

($z_0\in E$, $R$ is a Gaussian covariance operator, $M$ is a L\'evy
(spectral)

measure and $B$ is the unit ball in $E$);

$\nu=[z_0,R,M]$ means $\nu \in ID$ with the triple from its
L\'evy-Khintchine

formula;

$\nu^{\ast c}=[z_0, R, M]^{\ast c}:=[c \cdot z_0, c\cdot R, c\cdot
M]$ for $c>0$;

$(f\mu)$ denotes the image of a measure $\mu$ under a measurable
mapping $f$;

$(f\otimes g)(s_1,s_2):=f(s_1)\cdot g(s_2)$ (tensor product) for
$(s_1, s_2)\in S \times S$

and $f,g :S \to \Rset$;

$I^{h,r}_{(a,b]}$ is \emph{the random integral mapping} with \emph{a
space transform} function $h$,

a deterministic \emph{ monotone time change} $r$(an inner clock) and
the time

interval $(a,b]$; cf. (1) below;

$\mathcal{D}^{h,r}_{(a,b]}$ \emph{ is the domain} of the integral
mapping $I^{h,r}_{(a,b]}$;

$\mathcal{R}^{h,r}_{(a,b]}:=\,I^{h,r}_{(a,b]}(\mathcal{D}^{h,r}_{(a,b]})\subset
ID$ denotes the range of the mapping $I^{h,r}_{(a,b]}$;

$I^{h,r}_{(a,\infty)}(\nu)$ means weak limit of
$I^{h,r}_{(a,b]}(\nu)$ as $b\to \infty$;

$\Phi^{h,r}_{(a,b]}(y):= \log{\widehat{(I^{h,r}_{(a,b]}(\nu))}}(y)$
when $\Phi=\log\hat{\nu}(y)$

$I^{h,r}_{(a,b]}([z,R,M])=:[ z^{h,r}_{(a,b]}, R^{h,r}_{(a,b]},
M^{(h,r}_{(a,b]}]$, cf. (10);

\medskip
\medskip
\emph{0.2. \underline{Summary of results}}

\medskip
In section 1 the random integral and its basic properties are given.
In Theorem 1, in section 2, we proved that (some) mappings $I^{h,r}_{(a,b]}$
are isomorphisms, of the corresponding measure convolution subsemigroups of 
the semigroup $ID(E)$,
but not always (Remark 3).  An
alternative approach on $\Rset^d$, for retrieving the measures 
from random integral mappings, 
is given in Proposition 2. Then, in section 3, we discuss the domains
$\mathcal{D}^{h,r}_{(a,b]}$ of mappings $I^{h,r}_{(a,b]}$
(Propositions 3 - 6). In section 4, Theorem 2 shows that all
compositions of $I^{h,r}_{(a,b]}$ (iterated integral mappings) can
be expressed as a single  integral random mapping. Here the language
of tensor products and the notion of image measures are very
convenient. In section 5, the inverse mappings to $I^{h,r}_{(a,b]}$
are discussed in Theorem 3. However, they are no longer of the  
random integral form
$I^{h,r}_{(a,b]}$. Section 6, in particular Proposition 7, is
devoted to fixed points of  mappings $I^{h,r}_{(a,b]}$ and to the
role of the stable distributions. In section 7, the factorization
property of measures is discussed (Proposition 8). As a consequence
we get that the selfdecomposable (in other words class L
distributions) measures have the factorization property (Corollary 12). In the last
section (section 8) we illustrate our results on some new or
previously studied integral mappings and semigroups.

\medskip
\medskip
\textbf{1. A path-wise  random integral mappings.}

\medskip
\emph{1.1. \underline{Integrals on finite intervals}}

For an interval ${(a,b]}$ in a positive half-line, a
real-valued continuous of bound variation function $h$ on $(a,b]$, a
positive non-decreasing right-continuous (or non-increasing left-continuous) 
time change function $r$ on $(a,b]$ and a cadlag L\'evy stochastic processes $(Y_{\nu}(t), 0 \le
t < \infty)$, let us define via a formal integration by parts
formula the following \emph{random integral}
\begin{multline}
\int_{(a,b]} h(t)dY_{\nu}(r(t)):=
\\ h(b)Y_{\nu}(r(b))-h(a)Y_{\nu}(r(a))- \int_{(a,b]}Y_{\nu}(r(t)-)dh(t) \in E, 
\end{multline}
and the corresponding \emph{random integral mapping}
\begin{equation}
\qquad \qquad  \nu\to  I^{h,r}_{(a,b]}(\nu):=\mathcal{L}\big( \int_{(a,b]}
h(t)dY_{\nu}(r(t))\big) \in ID,  
\end{equation} 
with $\nu$ in its \emph{domain} $\mathcal{D}^{h,r}_{(a,b]}$ being a
subset of the class ID consisting of those measures $\nu$ for which
the integral (1) is well defined.
In that case,  the law in (1)  is  infinitely divisible; 
cf. Jurek-Vervaat (1983), Lemma 1.1.

\begin{rem}
(i) Note that if the random integral $\int_{(a,b]} h(t)dY_{\nu}(r(t))$ is well defined then so are
random integrals $\int_{(c,d]} h(t)dY_{\nu}(r(t)$, where $a\le c<d\le b$.

(ii) Since the random integral  $\int_{(a,b]} h(t)dY_{\nu}(r(t))$ is a functional of the   process on $(a,b]$, thus if  two L\'evy processes
$\bar{Y}_{\nu}$ and $Y_{\nu}$ have the same probability distribution,  that  is, $(\bar{Y}_{\nu}(t): t\ge
0)\stackrel{d}{=} (Y_{\nu}(t): t\ge 0)$, then
\begin{equation*}
 \ \ \ \  \ \ \qquad \mathcal{L}\big( \int_{(a,b]}
h(t)dY_{\nu}(r(t)\big)=\mathcal{L}\big( \int_{(a,b]}
h(t)d\bar{Y}_{\nu}(r(t)\big). \qquad \qquad \qquad   \qquad
\end{equation*}

(iii) Since L\'evy processes are semi-martingales the random integral (1)
can be defined as an Ito stochastic integral. However, for our purposes we do not need that generality of stochastic calculus.
\end{rem}

\medskip
\emph{1.2. \underline{L\'evy exponents (characteristic functions) of random mappings}}

\medskip
If $\nu \in \mathcal{D}^{h,r}_{(a,b]}$ and $I^{h,r}_{(a,b]}(\nu)$
have the L\'evy exponents $\Phi$ and $\Phi^{h,r}_{(a,b]}$,
respectively then, from already mentioned in Lemma 1.1 in Jurek-Vervaat(1983), we get
\begin{equation}
\Phi^{h,r}_{(a,b]}(y)=\int_{(a,b]} \Phi(h(t)y)dr(t),  \ \ y\in
E^{\prime} \ \ \ \mbox{(for non-decreasing $r$).}
\end{equation}
Similarly we have that
\begin{equation}
\Phi^{h,r}_{(a,b]}(y)= \int_{(a,b]} \Phi(-h(t)y)|dr(t)|, \ \  y\in
E^{\prime} \ \ \ \mbox{(for non-increasing $r$)},
\end{equation}
because for $0<u<w$, we have $\mathcal{L}(Y_{\nu}(u)-Y_{\nu}(w))=
(\nu^{-})^{\ast(w-u)}$ where $\nu^{-}:=\mathcal{L}(-Y_{\nu}(1))$. In
other words, $(-Y_{\nu}(t), t\ge0)\stackrel{d}{=}(Y_{\nu^-}(t),
t\ge0)$.

\medskip
\emph{1.3. \underline{Improper random integrals}}

\medskip
Integrals over intervals (a,b) or (a,$\infty$) or [a,b] and
others are defined as week limits of integrals over
intervals (a,b] in (1).  
Thus, the random integral $\int_{(a,\infty)} h(t)dY_{\nu}(r(t))$ is well-defined if and only if 
the function 
\[
E^{\prime}\ni y\to \int_{(a,\infty)} \Phi(h(t)y)dr(t)\in \Cset  \ \ \mbox{is a L\'evy exponent}. 
\]
\noindent Or equivalently, the three parameters $ z^{h,r}_{(a,\infty)}, R^{h,r}_{(a,\infty)}$ and
$M^{h,r}_{(a,\infty)}$ in L\'evy-Khintchine formula are well-defined; cf. (12) below.

\medskip
\emph{1.4. \underline{Different graphic notations}}

\medskip
Note that we have
\[
\qquad I^{h,\,
r}_{(a,b]}\,(\nu)=\mathcal{L}\Big(\int_0^{\infty}\,1_{(a,b]}\,(r^*(s))
h(r^*(s))\,dY_{\nu}(s)\Big) =
I^{\,\tilde{h}(s),\,\,s\,}_{(0,\infty)}(\nu) \equiv
I^{\tilde{h}}(\nu),  
\]
where $\tilde{h}(s):=1_{(a,b]}\,(r^*(s)) h(r^*(s))$ and $r^*$ is the
inverse function of $r$.

However, instead of that above graphicaly simpler notation, for a
greater flexibility of our considerations, we will keep the three
parameters: the time interval $(a,b]$, the space-valued normalization
$h$ and the inner time change $r$ in symbols and notions related to
the random integral mappings (1).

For improper random integrals with decreasing $r$ with
$r(a+)<\infty$ we have
$I^{h,\, r}_{(a,b]}\,(\nu)=I^{-h,\, r(a+)-r}_{(a,b]}\,(\nu)=I^{h,\,
r(a+)-r}_{(a,b]}\,(\nu^-)$,
that is,
\[
\int_{(a,b]} h(t)dY_{\nu}(r(t))\stackrel{d}{=}\int_{(a,b]}
h(t)dY_{\nu^-}(r(a+)-r(t)),
\]
and $t\to r(a+)-r(t)$ is a positive increasing function.

\medskip
\medskip
\textbf{2. Properties of  random integral mappings.}

\medskip
\emph{2.1. \underline{The isomorphism property}}

\medskip
\begin{thm}
Assume that $h(a):=h(a+)$, $r(a):=r(a+)$ exists in $\Rset$,  \qquad  $r$ is continuously differentiable  in $(a,b)$ and $h\neq 0$ on $[a,b]$. Then the mapping
\begin{equation}
 \mathcal{D}^{h,r}_{(a,b]}\ni\nu\to I^{h,r}_{(a,b]}(\nu) \,\in
\mathcal{R}^{h,r}_{(a,b]}
\end{equation}
is a continuous isomorphism between the corresponding measure convolution
semigroups. 
\end{thm}

\emph{Proof.}
Let $\rho:=I^{h,r}_{(a,b]}(\nu)$ and let  $\Phi^{h,r}_{(a,b]}$ and $\Phi$ be the L\'evy exponents of $\rho$ and $\nu$, respectively.
Note that for fixed $y \in E^{\prime}$, the function $[a,b]\ni s\to\Phi(h(s)y)r^{\prime}(s)$ is continuous.  By the Mean Value Theorem there exist $ c_y \in (a,b)$ such that
\begin{multline*}
\Im \Phi^{h,r}_{(a,b]}(y)=\int_{(a,b]}\Im\Phi(h(s)y)dr(s)\\ =\int_{[a,b]}\Im\Phi(h(s)y)r^{\prime}(s)ds=(b-a)\Im\Phi(h(c_y)y)r^{\prime}(c_y),
\end{multline*}
where $\Im$ stands for the imaginary part of a complex number.
Therefore we have that  $\Im \Phi(y)= [(b-a)r^{\prime}(c_y)]^{-1}\Im\Phi^{h,r}_{(a,b]}(y/h(c_y))$. Analogous equality holds for the real parts of $\Re\Phi$ and $\Re\Phi^{h,r}_{(a,b]}$. Consequently, we get that $\rho$ uniquely determines $\nu$, which proves the one-to-one property.

The homomorphism property of $I^{h,r}_{(a,b]}$, that is, the equality
\[
I^{h,\, r}_{(a,b]}\,(\nu_1\ast\nu_2))=I^{h,\,
r}_{(a,b]}\,(\nu_1)\ast I^{h,\, r}_{(a,b]}\,(\nu_2),
\]
in terms of the corresponding L\'evy exponents, follows from (3) or (4).

For the continuity, let us note that $0\leq |r(b)-r(a)|<\infty$  and
the cadlag property imply that functions $t\to Y(r(t))$ are bounded
and with at most countable many discontinuities; cf. Billingsley
(1968), Chapter 3, Lemma 1. Furthermore, the mapping
\begin{multline} D_E[a,b] \ni y\to
\int_{(a,b]}h(t)dy(r(t)):=
\\ h(b)y(r(b))-h(a)y(r(a))- \int_{(a,b]}y(r(t)-)dh(t) \in E,
\end{multline}
is continuous in Skorohod topology (for details see Billingsley
(1968), p. 121.). Furthermore, if $\nu_n \Rightarrow \nu$ then
$(Y_{\nu_{n}}(t), a \le t \le b) \Rightarrow (Y_{\nu}(t), 0\le t\le
b)$ in $D_E[a,b]$.  Consequently, we have
\[
\mathcal{L}\Big(\int_{(a,b]}h(t)dY_{\nu_{n}}(r(t))\Big)\Rightarrow\mathcal{L}\Big(\int_{(a,b]}h(t)dY_{\nu}(r(t))
\]
which proves the continuity of $I^{h,r}_{(a,b]}$ and completes the proof of Theorem 1.

\begin{rem} (i) In some specific cases  as
$I^{e^{-t},t}_{(0,\infty)}$ or $I^{t, t^{\beta}}_{(0,1]}$,  the
one-to-one property can be to proved by  Fourier or Laplace
transforms; cf. Jurek-Vervaat (1983) or Jurek (1988), (2007).

(ii) Weak convergence  continuity of the mappings $\nu\to I^{h,r}_{(a,b]}(\nu)$, for measures on  finite dimensional linear spaces, easily follows by the characteristic functional argument .
\end{rem}
\begin{cor}
For any $s>0$,\ $\nu \in
\mathcal{D}^{h,r}_{(a,b]}$ if and only if $\nu^{\ast
s}\in\mathcal{D}^{h,r}_{(a,b]}$,  and
\begin{equation*} 
I^{h,r}_{(a,b]}(\nu^{\ast s})= (I^{h,r}_{(a,b]}(\nu))^{\ast s} = (I^{h,\,sr}_{(a,b]}(\nu)).
\end{equation*}
For $u \in\Rset$ and the dilation operator $T_u$, $\nu \in
\mathcal{D}^{h,r}_{(a,b]}$ if and only if $T_u \nu \in
\mathcal{D}^{h,r}_{(a,b]}$, and
\begin{equation*}
T_u\big(I^{h,r}_{(a,b]}(\nu)  \big) =  I^{h,r}_{(a,b]}(T_u\nu)=I^{u
h,r}_{(a,b]}(\nu).
\end{equation*}
For bounded linear operator A on E and $\nu \in
\mathcal{D}^{h,r}_{(a,b]}$ we have that $A\nu\in\mathcal{D}^{h,r}_{(a,b]}$ and  
\[
A(I^{h,r}_{(a,b]}(\nu))=I^{h,r}_{(a,b]}(A\nu).
\]
\end{cor}
These are consequences of the formula (3) and (4).

\medskip
\emph{2.2. \underline{Convolution factors}}

\medskip
We say that probability measures $\mu$ on $E$ is \emph{a convolution
factor} of a measure $\rho$ if there exists a measure $\nu$ such
that $\mu\ast \nu=\rho$; in symbols we write $\mu\prec \rho$.
\begin{prop}
For $\nu \in \mathcal{D}^{h,r} _{(a,b]}$, the family $\{I^{h,r}
_{(a,x]}(\nu): a<x<b\}$ is sequentially shift convergent for
$x\uparrow c \le b$ or $x\downarrow c\ge a$.
\end{prop}
\emph{Proof}. Note that if $a<x_1<x_2< ... <x_n<... \uparrow c\le b$
then
\[I^{h,r} _{(a,x_1]}(\nu)\prec I^{h,r}_{(a,x_2]}(\nu)\prec ...
\prec I^{h,r} _{(a,b]}(\nu),
\]
and by Theorem 5.3 in Parthasarathy (1968) there exist sequence
$\delta_{z_n}$ and a measure $\rho$ such that
$I^{h,r}_{(a,x_n]}(\nu)\ast \delta_{z_n} \Rightarrow \rho$ as $n\to
\infty$. Similarly we argue in the remaining case.

\medskip
\medskip
\emph{2.3. \underline{Retriving the measure $\nu$}}

\medskip
Knowing  integrals (1) over  a family of intervals $(c,x]$, with
$x\downarrow c$, we can retrieve the measure $\nu$, as it  is seen
below.

\begin{prop} 
Let  $\nu\in ID(R^d)$, $r$ is continuously differentiable  and  there exists
$c\in(a,b)$ such that $h(c)\neq 0$ and $r^{\prime}(c)\neq 0$. Then
\begin{equation}
\mathcal{L}\Big(\int_c^x
\frac{h(t)}{h(c)}\,dY_{\nu}(\frac{r(t)}{r^{\prime}(c)})\Big)^{\ast\frac{1}{x-c}}\Rightarrow
\nu,         \ \ \mbox{as} \ x \downarrow c.
\end{equation}
\end{prop}

\emph{Proof.} The weak convergence in (6), for measures $\nu$ on $\Rset^d$, in
terms of L\'evy exponents is equivalent to the following claim
\[
\lim_{x\to c}\frac{1}{x-c}\,\int_c^x\,
\Phi\big(\frac{h(t)}{h(c)}\,y\big)\,
\frac{r^{\prime}(t)}{r^{\prime}(c)}dt=\Phi(y) \ \ \mbox{for all} \
y,
\]
that is obviously  true because of de l'Hospital rule.

\begin{cor}
For a measure $\nu\in ID(\Rset^d)$,
\[
\mathcal{L}\Big(\int_a^{a+1/n} \,h(t)\,dY_{\nu}(r(t))\Big)^{\ast\,
n}\Rightarrow T_{h(a+)}\nu^{\ast\,r^{\prime}(a+)} \ \mbox{as} \ n\to
\infty.
\]
\end{cor}
Note that $\nu$ is uniquely determined whenever $h(a+)\neq 0$ and $r^{\prime}(a+)\neq 0$.

\medskip
\medskip
\textbf{3. Domains $\mathcal{D}^{h,r}_{(a,b]}$ of integral mappings
$I^{h,r}_{(a,b]}$.}

\medskip
\emph{3.1. \underline{Domains on Banach spaces $E$.}}

\begin{prop} In order that $\mathcal{D}^{h,r}_{(a,b]}=ID(E)$ it is
necessary and sufficient that integrals
$\int_{(a,b]}y(r(t)-)\,dh(t)$ exists for all $y\in D_E[a,b]$.
In particular, if
$|r(b)-r(a+)|<\infty$ then $\mathcal{D}^{h,r}_{(a,b]}=ID(E)$.
\end{prop}
\emph{Proof.} Because L\'evy processes $Y$ are cadlag and the random
integrals (1) are defined by the formal integration by parts, we infer the claim 
concerning the first part.

Since the range of $r$ is bounded then using the fact that  cadlag
functions, on bounded intervals,  are integrable  (cf. Billingsley
(1968), p.121) we get that the integral in (1) is well-defined.
This concludes the second claim.

\medskip
In terms of the  parameters in L\'evy-Khintchine representation, domains of random integrals  are  characterized as follows:
\begin{prop}
A measure $\nu=[z,R,M]$ is in the domain $\mathcal{D}^{h,r}_{(a,b]}$ if and only if the following holds
\begin{equation}
\int_{(a,b]}|h(t)||dr(t)|<\infty, \ \mbox{if} \ z \neq 0; \
\int_{(a,b]}h^2(t)|dr(t)|<\infty,  \ \mbox{if} \ R \neq 0,
\end{equation}
and for the $\sigma$-algebra $\mathcal{B}_0$ of Borel subsets of
$E\setminus{\{0\}}$, the set function
\begin{equation}
\mathcal{B}_0 \ni A\to \int_{(a, b]}[T_{h(t)}M(A)]|dr(t)|\ \mbox{is
a L\'evy spectral measure on E}.
\end{equation}
Moreover, if $I_{(a,b]}^{h,r}(\nu)$ is determined by the triple
$[z^{h,r}_{(a,b]}, R^{h,r}_{(a,b]}, M^{(h,r}_{(a,b]}]$ and  $r$ is
nondecreasing then
\begin{multline}
(i) \ z^{h,r}_{(a,b]}=(\int_{(a,b]}h(t)dr(t))\cdot z +
\int_{(a,b]}h(t) \int_{E\setminus\{0\}}[1_{B}(h(t)x)-1_{B}(x)] x
M(dx) dr(t); \\
(ii) \ \ R^{h,r}_{(a,b]}=(\int_{(a,b]}h(t)^2 dr(t))\cdot R;\\
(iii)\ M^{h,r}_{(a, b]}(A)= \int_{(a, b]}[T_{h(t)}M(A)]dr(t)
=\int_{(a, b]}\int_{E\setminus \{0\}} 1_A(h(t)x)\,M(dx)\,dr(t),
\end{multline}
\end{prop}
\emph{Proof.}  From formulas in Section 0.1, (2) and from the
uniqueness of the triple (shift vector, Gaussian covariance and
L\'evy spectral measure), in the L\'evy-Khintchine formula, we get the
above claims and the three formulas in (9).
\begin{cor}
If $ M^{h,r}_{(a, b]}$ is a L\'evy spectral measure on $E$ then
\[
\int_{(a,b]}(1\wedge h^2(s))|dr(s)|<\infty.
\]
\end{cor}
\emph{Proof.} For $y\in E^{\prime}$ and the mapping $\pi_y (x) :
=<y,x>$ ( $x\in E$), the image measure $\pi_y ( M^{h,r}_{(a, b]})$
is a L\'evy spectral measure on $\Rset$. Since for positive s and t
we have $(1\wedge s)(1\wedge t)\le (1\wedge st)$ therefore we have
that
\begin{multline*}
(\int_{(a,b]}(1\wedge h^2(s))|dr(s)|)\cdot(\int_E(1\wedge
<y,x>^2)\,dM(x)) \\ \le \int_{(a,b]}\int_E (1\wedge
<y,h(s)x>^2)\,dM(x)|dr(s)|\\ = \int_E (1\wedge
<y,u>^2)\,M^{h,r}_{(a, b]}(du) = \int_\Rset (1\wedge w^2)\, (\pi_y
M^{h,r}_{(a, b]})(dw)<\infty,
\end{multline*}
which concludes the proof.

\begin{cor}
Let  h be any real-valued function on $(a,b]$,  let  r be monotone on the interval (a,b] 
such that $|r(b)-r(a+)| <\infty$ and  
let L\'evy (spectral)  measures M and N be  such that $[0,0,M]$ and $[0,0,N]$ 
are in the domain $\mathcal{D}^{h,r}_{(a,b]}$. 
Then  $M^{h,r}_{(a, b]}=N^{h,r}_{(a, b]}$ implies that M=N.
\end{cor}
\emph{Proof.} Assume that $r$ is non-decreasing and $r(b)-r(a+)<\infty$. 
Then for any Borel subset $A \subset E\setminus\{0\}$ and bounded away from zero,
i.e., $0\notin \bar{A}$ (closure),
\[
 M^{h,r}_{(a, b]}(A) -N^{h,r}_{(a, b]}(A)  =
\int_{E\setminus\{0\}}[\int_{(a,b]}1_A(h(t)x)dr(t)](M-N)(dx)=0,
\]
which can be extended (note that $r(b)-r(a+)<\infty$) to the following
\[
\int_{E\setminus\{0\}}[\int_{(a,b]} g_0(h(t)x)dr(t)](M-N)(dx)=0, \ 
\mbox{for all} \  \ g_0\in C_{b\,0}^+(E),
\]
where $C_{b\,0}^+(E)$ stands
for the family of all  functions that are positive, continuous bounded and vanishing in 
some neighbourhoods of zero. Since the expression in the square bracket is always positive
we conclude that M=N, which completes the proof.

\begin{rem}
Choose two different L\'evy spectral measures $M$ i $N$ concentrated on $(0,\infty)$ such that
\[
\int_{(0,\infty)}x^2M(dx)=\int_{(0,\infty)}x^2N(dx)<\infty,
\]
and the time change function $r(t)= t^{-2}$.
Then for each $v>0$,
\begin{multline*}
M^{t,\,\,t^{-2}}_{(0,\infty)}(x>v)= 2
\int_0^{\infty}\int_{(0,\infty)}1_{(x>v)}(tx)t^{-3}dtM(dx) \qquad (w:=tx)\\
= \int_{(0,\infty)}x^2M(dx)\,\, 2\int_v^{\infty}w^{-3}dw =
\int_{(0,\infty)}x^2M(dx)\,v^{-2}\\ 
=\int_{(0,\infty)}x^2N(dx)\,v^{-2}=N^{t,\,\,t^{-2}}_{(0,\infty)}(x>v).
\end{multline*}
So, for different L\'evy measures $M$ and $N$ we got equality 
$M^{t,\,\,t^{-2}}_{(0,\infty)}=N^{t,\,\,t^{-2}}_{(0,\infty)}$. However, the functions  $h(t)=t$ and $r(t)=t^{-2}$
do not satisfy the inegrability condition from Corollary 3. 
Thus $ M^{t,\,\,t^{-2}}_{(0,\infty)}$ is not a L\'evy (spectral)   measure.
\end{rem}

Here are some sufficient conditions (for symmetrized measures $\gamma^{\circ}$) to be in a domain
of a random integral mapping.
\begin{prop}
(i) For $ 0<p\le2$, if the integral $\int_{(a,b]}|h(t)|^p |dr(t)|$
exists then  all symmetric p-stable measures $\gamma^{\circ}_p$ on
$E$ are in  the domain $\mathcal{D}^{h,r}_{(a,b]}$.

\noindent (ii) If a positive Borel measure $N$ on $E$  integrates
the function $||x||$ then $N$ is a L\'evy spectral measure and
$\nu=[z,0,N]$ has finite first moment. Moreover, if $
\int_{(a,b]}|h(t)|\,|dr(t)|<\infty$ then
$[z,0,N]\in\mathcal{D}^{h,r}_{(a,b]}$.
\end{prop}

\emph{Proof.} For (i), recall that L\'evy exponents of symmetric
p-stable non Gaussian measures for $0<p<2$ are of the form
\[
\Phi(y)\equiv - \log \hat{\gamma^{\circ}_p}
(y)=\int_{\{||x||=1\}}|<y,x>|^p\,m(dx)
\]
for some finite measure m on the unit sphere; cf. Araujo and Gin\'e
(1980), Chapter III, Theorem 6.16. Hence and from (3)
\[
y\to \int_{(a,b]}\Phi((h(t)y) dr(t)= \int_{(a,b]}|h(t)|^p\,dr(t)
\int_{\{||x||=1\}}|<y,x>|^p\,m(dx)
\]
is also a L\'evy exponent and thus the random integral is well
defined. The case of  symmetric Gaussian ($p=2$) follows from
Corollary 4 (ii).

For part (ii), since $\int_{E\setminus \{0\}} (1\wedge
||x||)\,N(dx)<\infty$, therefore  $N$ is a L\'evy measure by
Araujo-Gine (1980), Chapter 3, Theorem 6.3. Since also

\noindent $\int_{(||x||>1)}\,||x||N(dx)<\infty$ we conclude that
$\nu$ has finite first moment. Furthermore for measure
$N^{h,r}_{(a,b]}$ given by (14) we have
\begin{multline}
\int_{E\setminus \{0\}} (1\wedge
||x||)\,N^{h,r}_{(a,b]}(dx)=\int_{(a,b]}\int_{E\setminus \{0\}}
(1\wedge
|h(t)||x||)N(dx)|dr(t)|\\
 \le (\int_{(a,b]}|h(t)||dr(t)|)\ (\int_E ||x||N(dx))<\infty,
\end{multline}
and again by Theorem 6.3 in  Chapter 3 in Araujo-Gin\'e (1980) we
conclude that $N^{h,r}_{(a,b]}$ is a L\'evy spectral measure. Thus
$\nu\in\mathcal{D}^{h,r}_{(a,b]}$ and the proof is completed.

\medskip
\emph{3.2.\underline{Domains on Hilbert space $H$}}

\medskip
On real separable Hilbert spaces we have complete characterization
of covariance operators and more importantly, for the considerations
here, we know that
\[
\mbox{M is a L\'evy measure on H} \ \mbox{iff}  \ M\{0\}=0 \
\mbox{and}\ \int_H\, (1\wedge ||x||^2)\,M(dx)<\infty,
\]
cf. Parthasarathy (1968), Chapter VI.  With the above and
Proposition 4 we have
\begin{cor}
In order that  a mesure $\nu =[z, R,M]\in \mathcal{D}^{h,r}_{(a,b]}(H)$ it is necessary and 
sufficient that \\
$(i)\ \ \int_{(a,b]}|h(t)|\,|dr(t)|<\infty, \ \
\mbox{provided} \ z \neq 0,
\\ (ii) \ \ \int_{(a,b]}h^2(t)|dr(t)|<\infty , \  \  \mbox{provided} \ R
\neq 0, \\
(iii) \ \  \int_{(a, b]}\int_{E\setminus \{0\}} (1\wedge
|h(t)|^2||x||^2)\,M(dx)\,|dr(t)|<\infty,  \ \ \mbox{provided} \  M
\neq 0.$
\end{cor}

\begin{rem} Since for all positive s and t,
$(1\wedge s)(1\wedge t)\le (1\wedge st)$, from the above condition
(iii) we infer that  if $M^{h,r}_{(a,b]}$ is a L\'evy spectral
measure (on H) then so is $M$ and
\[
\int_{(a, b]} (1\wedge |h(t)|^2)\,|dr(t)|<\infty.
\]
Therefore, if $[0,0,M]\in\mathcal{D}^{h,r}_{(a,b]}(H)$  then it is necessary that  the triple: an interval $(a,b]$
and functions $h$ and $r$, satisfies the above integrability condition.
\end{rem}

\begin{prop}
For  triples $(a,b], h$ and $r$  satisfying the conditions (i) and
(ii) from Corollary 4, all infinitely divisible measures with finite
second moment  are in their domains, that is, $ID_2(H)\subset
\mathcal{D}^{h,r}_{(a,b]}$, for arbitrary Hilbert space H.
\end{prop}
\emph{Proof.} In view of Jurek-Smalara (1981) or Proposition 1.18.13
in Jurek-Mason(1993) or Theorem 25.3 in Sato (1999) we know that
$\nu=[z,R,M]\in ID_2(H)$ if and only if
$\int_{(||x||>1)}||x||^2M(dx)<\infty$. Since
\[
\int_{(a, b]}\int_H (1\wedge |h(t)|^2||x||^2)\,M(dx)\,|dr(t)\le
\int_{(a,b]}h^2(t)|dr(t)|\, \ \int_H||x||^2M(dx)< \infty,
\]
( on H, L\'evy measure M always integrates $||x||^2$ in the unit
ball), we conclude that $\nu\in \mathcal{D}^{h,r}_{(a,b]}$, which
completes the proof.

\medskip
\medskip
\textbf{4. Compositions of random integral mappings
$I^{h,r}_{(a,b]}$}

\medskip
\emph{4.1. \underline{Equivalent mappings}}

\medskip
We say that two integral mappings $I^{h,r}_{(a,b]}$ and
$I^{h_1,r_1}_{(a_1,b_1]}$  are \emph{equivalent} if
\begin{equation}
\mathcal{D}^{h,r}_{(a,b]}=\mathcal{D}^{h_1,r_1}_{(a_1,b_1]} \ \ \
\mbox{and} \ \  \ I^{h,r}_{(a,b]}(\mathcal{D}^{h,r}_{(a,b]})=
I^{h_1,r_1}_{(a_1,b_1]}(\mathcal{D}^{h_1,r_1}_{(a_1,b_1]}),
\end{equation}
and we write $I^{h,r}_{(a,b]} = I^{h_1,r_1}_{(a_1,b_1]}$. In terms
of L\'evy exponents the above means that
\[
\int_{(a_1,b_1]}\Phi(h_1(t)y)dr_1(t)=\int_{(a_2,b_2]}\Phi(h_2(t)y)dr_2(t),
\ \mbox{for all} \ y\in E^{\prime}
\]
for L\'evy exponents $\Phi$  (measures) in appropriate domains.

\begin{rem}
Mappings $I^{e^{-t},t}_{(0,\infty)}$ and $I^{s, \,-\log s}_{(0,1)}$
are equivalent. Similarly, $I^{t,t^{\be}}_{(0,1]}$ and
$I^{t^{1/{\be}},t}_{(0,1]}$, for $\be>0$. This follows from above
without specifying the domains.
\end{rem}

\medskip
\emph{4.2. \underline{Iterated random integral mappings}}

\medskip
Below let the time change $r(t),\, a<t\le b,$ be either $\rho\{s:
s>t\}$ or $\rho\{s: s\le t\}$ for some positive, possibly infinite,
measure $\rho$ on a positive half-line.

\medskip
For functions $h_1,...,h_m$, intervals $(a_1, b_1],..., (a_m,b_m]$
and measures $\rho_1,...,\rho_m$ let us define
\begin{multline}
\textbf{h}:= h_1 \otimes...\otimes h_m,\ \ \mbox{(tensor product of
functions)} \\ \mbox{i.e.} \ \
\textbf{h}(t_1,t_2,...,t_m):=h_1(t_1)\cdot
h_2(t_2)\cdot...\cdot h_m(t_m), \ \mbox{where} \ a_i <t_i  \le b_i; \\
\textbf{(a,b]}:=(a_1,b_1]\times ...\times (a_m,b_m], \
\boldsymbol{\rho}:= \rho_1\times...\times \rho_m  \ \mbox{(product
measure)}
\end{multline}

\begin{thm}
Let functions $h_i$, measures $\rho_i$ (given by increments of
functions $r_i$) and intervals $(a_i,b_i]$, for $i=1,2,...,m$, be as
above.

If the image $\textbf{h((a,b])}=(c,d]\subset \Rset^+$ and $\nu\in
ID(E)$ is from an appropriate domain then we have
\begin{equation}
I_{(a_1, b_1]}^{h_1, \rho_1}(I_{(a_2, b_2]}^{h_2,
\rho_2}(...(I_{(a_m, b_m]}^{h_m,\rho_m}(\nu))))= I_{(c,d]}^{t, \,
\textbf{h}\,\,\boldsymbol{\rho}}(\nu)
\end{equation}
where $\textbf{h}\ \boldsymbol{\rho}$ is the image of the product
measure $\boldsymbol{\rho}=\rho_1\times ... \times \rho_m$ under the
mapping $\textbf{h}:=h_1 \otimes...\otimes h_m$.
\end{thm}
\emph{Proof.} For $\nu\in\mathcal{D}^{h,r}_{(a,b]}$ and its L\'evy
exponent $\Phi$ let us define the (script) mapping
$\mathcal{I}^{h,r}_{(a,b]}$ as follows
\begin{equation}
\mathcal{I}^{h,r}_{(a,b]}(\Phi)(y):= \Phi^{h,r}_{(a,b]}=
\int_{(a,b]}\Phi(\pm h(s)y)d(\pm)r(s),
\end{equation}
where the sign minus is in the case of decreasing time change $r$.
Then to justify (14) it is enough to notice that
\begin{multline}
\mathcal{I}_{(a_1, b_1]}^{h_1, \rho_1}(\mathcal{I}_{(a_2,
b_2]}^{h_2, \rho_2}(...(\mathcal{I}_{(a_m,
b_m]}^{h_m,\rho_m}(\Phi))))(y)\\ =
\int_{(a_1,b_1]}\int_{(a_1,b_2]}...\int_{(a_m,b_m]}\Phi\big(h_1(t_1)\,h_2(t_2)\,...\,h_m(t_m)\,y))\big)\,dr_m(t_m)...dr_2(t_2)dr_1(t_1)\\
=\int_{(\textbf{a},\textbf{b}]}\Phi\big(h_1 \otimes...\otimes
h_m(s)\,y\big) \boldsymbol{\rho}(ds) =
\int_{(c,d]}\Phi(t\,y)(\textbf{h} \boldsymbol{\rho})(dt),
\end{multline}
which follows from the Fubini and the image measure theorems.

In view of the definitions of the tensor product and the product
measures we have
\[
h_1 \otimes...\otimes
h_m\,(\rho_1\times\rho_2\times...\times\rho_m)=h_{\sigma(1)}
\otimes...\otimes
h_{\sigma(m)}\,(\rho_{\sigma(1)}\times\rho_{\sigma(2)}\times...\times\rho_{\sigma(m)})
\]
for any permutation $\sigma$ of $1,2,...m$. Hence
\begin{cor}
Random integrals $I^{h_i,\rho_i}_{(a_i,b_i]}, i=1,2,...,m$, commute
on the domain
$\mathcal{D}^{t,\,\textbf{h}\,\,\boldsymbol{\rho}}_{(c,d]}$, where
$\boldsymbol{\rho}=\rho_1\times...\times\rho_m$ and
$\textbf{h}:=h_1\otimes ...\otimes h_m$.
\end{cor}

In case of probability measures $\rho_i$, the time change function
$r$ is a cumulative probability distribution and we have
\begin{cor}
Let assume that $r_i(t):=\rho_i(\{s\in(a_i,b_i] : a<s\le t\})$,
where $\rho_i$ are probability measure on $(a_i,b_i]$ that are
distributions of random variables $Z_i$, for $1\le i \le m$. If
$Z_1, Z_2,..., Z_m$ are stochastically independent then
\[
r(t):=\textbf{h}\,\,\boldsymbol{\rho}(s\le t)=P[h_1(Z_1)\cdot ...
\cdot h_m(Z_m)\le t].
\]
\end{cor}
The above we can apply, for instance, to  $h_i(t):=|t|$ on positive
half-line and independent standard normal variable $Z_i$. That case
was investigated in $\Rset^d$ by Aoyama (2009) via polar
decomposition of L\'evy spectral measures.

\medskip
\emph{4.3. \underline{Inclusion of ranges of integral mappings}}

If a random mapping is a composition of other mappings we may infer
some inclusions of their ranges. Namely we have
\begin{cor} If  an equality $I^{h,r}_{(a,b]}=
I^{h_1,r_2}_{(a_1,b_1]}\circ I^{h_2,r_2}_{(a_2,b_2]}$ (a
composition) holds on the domain $\mathcal{D}^{h,r}_{(a,b]}$ then we
have
\[
\mathcal{R}^{h,r}_{(a,b]}\equiv
I^{h,r}_{(a,b]}(\mathcal{D}^{h,r}_{(a,b]})\subset
I^{h_1,r_1}_{(a_1,b_1]}(\mathcal{D}^{h_1,r_1}_{(a_1,b_1]})\cap
I^{h_2,r_2}_{(a_2,b_2]}(\mathcal{D}^{h_2,r_2}_{(a_2,b_2]})=\mathcal{R}^{h_1,r_1}_{(a_1,b_1]}\cap
\mathcal{R}^{h_2,r_2}_{(a_2,b_2]}
\]
\end{cor}
\emph{Proof.} From the equality of the above mappings we get
\[
I^{h,r}_{(a,b]}(\mathcal{D}^{h,r}_{(a,b]})=I^{h_1,r_1}_{(a_1,b_1]}\big(
I^{h_2,r_2}_{(a_2,b_2]}(\mathcal{D}^{h,r}_{(a,b]})\big) \ \mbox{ and
hence}\ \ I^{h_2,r_2}_{(a_2,b_2]}(\mathcal{D}^{h,r}_{(a,b]})\subset
\mathcal{D}^{h_1,r_1}_{(a_1,b_1]}.
\]
Therefore $I^{h,r}_{(a,b]}(\mathcal{D}^{h,r}_{(a,b]})\subset
I^{h_1,r_1}_{(a_1,b_1]}(\mathcal{D}^{h_1,r_1}_{(a_1,b_1]})$. Because
of the commutativity we get
$I^{h,r}_{(a,b]}(\mathcal{D}^{h,r}_{(a,b]})\subset
I^{h_2,r_2}_{(a_2,b_2]}(\mathcal{D}^{h_2,r_2}_{(a_2,b_2]})$, which
completes a proof.

\medskip
\medskip
\emph{4.4. \underline{An example of application of Theorem 2}}

\medskip
\begin{lem}
Let $ h_1(t):= e^{-t},  r_1(t):=t, h_2(s):=s$ and
$r_2(s):=1-e^{-s}$, \qquad  \ $0<s, t<\infty$. Then the corresponding measures are:
$d\rho_1(t)=dt, \, d\rho_2(s)=e^{-s}ds$ and
$d\boldsymbol{\rho}(t,s)= d(\rho_1\times
\rho_2)(t,s)=e^{-s}\,dt\,ds$. Finally, for the image measure
$\textbf{h}\,\,\boldsymbol{\rho}(dw)=(h_1\otimes h_2) (\rho_1\times
\rho_2)(dw)=\frac{e^{-w}}{w}dw$.
\end{lem}
\emph{Proof.} For Borel measurable, bounded and non-negative
functions $g$ we have
\begin{multline*}
\int_0^{\infty}g(u)(h_1\otimes h_2) (\rho_1\times
\rho_2)(du)=\int_0^{\infty}\int_0^{\infty}g((h_1\otimes
h_2)(t,s))\rho_1(dt)\rho_2(ds)\\
=\int_0^{\infty}\int_0^{\infty}g(e^{-t}\,s)dt\,e^{-s}ds=
\int_0^{\infty}(\int_0^{s}g(w)\frac{1}{w}\,dw)\,e^{-s}ds=
\int_0^{\infty}g(s)\,\frac{e^{-s}}{s}\,ds,
\end{multline*}
which completes the proof of Lemma 1.

From Theorem 2, Corollary 6, Lemma 1 and Jure-Vervaat (1983) we conclude that
\begin{cor} For $\nu\in ID_{\log}$  we have
\[
 I^{t,\,1-e^{-t}}_{(0,\infty)}\big(I^{e^{-s},\, s}_{(0,\infty)}(\nu)\big)=
I^{e^{-s},\,s}_{(0,\infty)}\big(I^{t,\,1-e^{-t}}_{(0,\infty)}(\nu)\big)=I^{-w,\,\Gamma(0;w)}_{(0,\infty)}(\nu)=I^{w,\,\Gamma(0;w)}_{(0,\infty)}(\nu^-)
\]
Moreover, $\Gamma(0; w)= (h_1\otimes h_2) (\rho_1\times
\rho_2)(\{x:x>w\})=\int_w^{\infty}\frac{e^{-s}}{s}\,ds$ for $w>0$.
\end{cor}

\begin{rem} (a) For the Euler constant \textbf{C} we have
\[
- \Gamma(0;w)=Ei(-w)= \textbf{C}+ ln\,w
+\int_0^w\frac{e^{-t}-1}{t}dt, \ \mbox{for} \  w>0,
\]
where $Ei$ is the special \emph{exponential-integral} function; cf.
Gradshteyn-Ryzhik (1994), formulae 8.211 and 8.212.

(b) Recall that the class 
$I^{t,\,1-e^{-t}}_{(0,\infty)}(ID)\equiv \mathcal{E}$ was introduced
in Jurek (2007), where the mapping $I^{t,\,1-e^{-t}}_{(0,\infty)}$
was denoted by $\mathcal{K}^{(e)}$; ((e) for exponential cumulative
distribution function). More importantly, the class $\mathcal{E}$
was related to the class of Voiculescu $\boxplus$ free-infinitely
divisible measures; cf. Corollary 6 in Jurek (2007). Note also that
$I^{t,\,1-e^{-t}}_{(0,\infty)}=I^{-\log s,\,s}_{(0,1]}$ and thus it
coincides with the upsilon mapping $\Upsilon$ studied in
Barndorff-Nielsen, Maejima and Sato (2006).

(c) Similarly $I^{e^{-s},\,s}_{(0,\infty)}(ID_{\log})\equiv L$
coincides with the L\'evy class of selfdecomposable probability
measures; cf. Jurek-Vervaat (1983), Theorem 3.2 or 
Jurek-Mason (1993), Theorem 3.6.6.

(d) Finally we get identity
$I^{e^{-s},\,s}_{(0,\infty)}\big(I^{t,\,1-e^{-t}}_{(0,\infty)}(ID_{\log})\big)
\equiv T$, which is the Thorin class; cf. Grigelionis (2007),
Maejima and Sato (2009) or Jurek (2011).
\end{rem}

\medskip
From Corollary 9 and Remark 5(d) we infer that
\begin{cor}
For the three classes: Thorin class T, L\'evy class L
(selfdecomposable measures) and $\mathcal{E}$ we have that $T\subset
L\cap\mathcal{E}$
\end{cor}
(This inclusion, on $\Rset^d$, was first noticed in
Barndorff-Nielsen, Maejima and Sato (2006) and also in Remark 2.3 in
Maejima-Sato (2009) but by using completely different methods.)

\medskip
\medskip
\textbf{5. The identity and the  inverse of a random integral mapping.}

\medskip
\emph{5.1. \underline{Zero  random integral mapping} $I^{h,r}_{(a,b]}$}

\medskip
If $h(t)\equiv 0$ or $r(t)\equiv r_0$ (constant function) then
\[
\int_{(a,b]} 0\, dY_{\nu}(r(t))=0=\int_{(a,b]} h(t)dY_{\nu}(r_0), \ \
\mbox{for all} \ \nu\in ID,
\]
are the zero mappings; $I^{0,r}_{(a,b]}(\nu)=I^{h,
r_0}_{(a,b]}(\nu)=\delta_0$. To exclude that trivial case we assume that
the three parameters: an interval
$(a,b]$ and functions $h, r$,  satisfy the basic condition
$0<\int_{(a,b]}|h(t)|\,|dr(t)|$.

\noindent On the other hand, the condition $\int_{(a,b]}|h(t)|\,|dr(t)|<\infty,$
guarantees that the
degenerate L\'evy process $Y(t):=ta$ (a fixed) can be used as
integrators in the integrals (1); cf. formula (8) in Proposition 4.

\medskip
\emph{5.2. \underline{Identity random integral mapping} $I^{h,r}_{(a,b]}$}

Note that whenever $0< r_0:=|r(b)-r(a+)|<\infty$ and $h(t)\equiv 1$
(constant) then
\begin{equation}
 I^{ 1,\, r(t)/r_0}_{(a,b]}(\nu))=\nu \ \ \mbox{for all} \
\nu\in ID.
\end{equation}
So all mappings $I^{1,\, r(t)/r_0}_{(a,b]}$ play a role of the
neutral element (identity mapping), under the composition,  in the family of all
integral mappings.  In fact, for $r(.)/r_0$ one may take any time
change whose increment over the interval $(a,b]$ is equal to 1.

Similarly, if $\delta_u(t):=1_{[u,\infty)}(t)$, $h(u)\neq0$ (u is
fixed) and $u \in (a,b)$ then from (2) or (3) we have
\begin{equation}
I^{h/h(u),\delta_u}_{(a,b]}(\nu)=\nu
\end{equation}
and thus they also play the role of the identity mapping.
\begin{rem}
In the case of (17),  the time change can be  any strictly monotone function
while the space change $h$ is trivial. In the  case of (18), the space change is quite
arbitrary but time change $r$ is one point jump function.  
\end{rem}
Integrals (18) and (19) are equivalent and will be called \emph{the identity mappings} in the space of all random integral mappings $I^{h,r}_{(a,b]}$.

\medskip
\medskip
\emph{5.3. \underline{The inverse of a random integral mapping.}}

\medskip
Under the conditions in Theorem 1, there exists the inverse  of the mapping
$I^{h,r}_{(a,b]}$ for which we have
\begin{thm}  If the mapping
$(I^{h,r}_{(a,b]})^{-1}\, : \, \mathcal{R}^{h,r}_{(a,b]}\equiv
(I^{h,r}_{(a,b]}(\mathcal{D}^{h,r}_{(a,b]}))\to
\mathcal{D}^{h,r}_{(a,b]}$  exits then it is an isomorphism between the
corresponding subsemigroups of ID. However, it is not of the
random integral mapping form, unless it is the identity mapping.
\end{thm}
\emph{Proof.} The isomorphism property of the inverse mapping is a
consequence of the fact that $I^{h,r}_{(a,b]}$ is an isomorphism by
Theorem 1.

Now suppose that the inverse of a non-trivial  mapping $I^{h,r}_{(a,b]}$ is,
indeed, of an integral form $I^{h_1,r_1}_{(a_1,b_1]}$. Then by
Theorem 2
\begin{equation*}
I^{h,r}_{(a,b]}(I^{h_1,r_1}_{(a_1,b_1]}(\nu))=I^{s,
r_2}_{(c,d]}(\nu)= \nu,
\end{equation*}
where
\begin{equation*}
dr_2(t) =d(h \otimes h_1)(dr\times dr_1)(t)=d\delta_u(t) \ \
\mbox{for some fixed} \ \ u\in (c,d],
\end{equation*}
 and $(c,d]=(h\otimes h_1) ((a,b]\times(a_1,b_1])$.  In other words, for all continuous functions $g$ on $[c,d]$ we have
\[
\int_{(a,b]}\int_{(a_1,b_1]}g(h(t)h_1(s))dr(t)dr_1(s)=g(u).
\]
Hence, either $h(t)\cdot h_1(s)= u$ (constant) for all $t\in (a,b]$
and $s\in (a_1,b_1]$  and $|r(b)-r(a+)| |r_1(b_1)-r_1(a_1+)|=1$ or
$dr\times dr_1=\delta_t \times \delta_s$ and $h(t)h_1(s)=u$.
Consequently, in the first case both $h$ and $h_1$ are constant that
contradicts the assumption that $I^{h,r}_{(a,b]}$ is non-trivial
mapping. Similarly, in the second case $r$ and $r_1$ are Dirac
measures and therefore $I^{h,r}_{(a,b]}$ is an identity mapping. Thus this completes a proof of Theorem 3.

\medskip
\textbf{6. Fixed points (eigenfunctions) of $I^{h,r}_{(a,b]}$}

\medskip
\emph{6.1. \underline{Definition of fixed points}}

We will say that an infinitely divisible measure $\rho$ is \emph{a
fixed point of an integral mapping $I^{h,r}_{(a,b]}$}, if
\begin{equation}
I^{h,r}_{(a,b]}(\rho)= \rho^{\ast c}\ast\delta_z, \ \qquad \mbox{for
some $c>0$ and $z\in E$.}
\end{equation}
Equivalently, in terms of L\'evy exponents, using (15)
\begin{equation}
\mathcal{I}^{h,r}_{(a,b]}(\Phi)(y)=c\,\Phi(y) +i<y,z> \ \mbox{for
all} \  y\in E^{\prime}.
\end{equation}

\begin{rem}
i)  Remark 5.2 in Jurek-Vervaat (1983) explains  why in the
definition (20) we have taken  $\nu^{\ast c}$ instead of the more natural
$T_c\nu$ (multiplying of a corresponding a random variable by  a
constant).

ii) Note that (20) reads that $\Phi$ is an eigenfunction of the
mapping $\mathcal{I}^{h,r} _{(a,b]}$ acting on the positive cone of
(symmetric)  L\'evy exponents, provided we ignore the shift part.
\end{rem}

\medskip
\emph{6.2. \underline{Stable measures}}

\medskip
Let us recall here one of the many equivalent definitions of stable
distributions. Namely, we will  say that \emph{$\gamma$ is a stable
probability measures} if there exists a parameter $0<p\le2$  such for each  
$t>0$ there exists a vector $z(t)\in E$ such that
\begin{equation}
t^p\,\log \hat{\gamma}(y)= \log \hat{\gamma}(ty)+i<y,z(t)> \ \
\mbox{for all}\ y \in E^{\prime};
\end{equation}
cf. Jurek (1983), Theorem 3.2. or Linde (1983) or Theorem 4.14.2 in Jurek-Mason (1993).
 We will write
$\gamma_p\equiv \gamma$ if the above holds and say that it is
\emph{a p-stable probability measure}. Furthermore, we say that $\gamma_p$ is
\emph{strictly stable}, if $z(t)\equiv 0$ in (21).

\medskip
\emph{6.3. \underline{Fixed points of the mapping $I^{h,r}_{(a,b]}$}}

\begin{prop}
In order that p-stable measure $\gamma_p$ be a fixed point of the
mapping   $I^{h,r}_{(a,b]}$ it is necessary and sufficient that
$0<\int_{(a,b]}|h(t)|^p |dr(t)|<\infty$ .
\end{prop}
\emph{Proof.} Because of the shift z in (20), if it enough to
consider only strictly stable measures. In that case, using (21), we have
\begin{multline}
\log\widehat{I^{h,r}_{(a,b]}(\gamma_p)}(y)\\ 
=\int_{(a,b]}\log\hat{\gamma_p}(h(t)y)\,dr(t)
=[\int_{(a,b]}\,|h(t)|^p\,dr(t )] \,\log\hat{\gamma_p}(y),
\end{multline}
that is, p-stable probability measures $\gamma_p$ are fixed points of the mapping
$I^{h,r}_{(a,b]}$ with the constant  $c:=\int_{(a,b]}|h(t)|^p |dr(t)|$, which
completes the proof.

\medskip
Let denote by $\mathcal{S}$ the set of all stable measures. Then we
get
\begin{cor}
For the class $\mathcal{S}$ of all stable measures
\[
\big[I^{h,r}_{(a,b]}(\mathcal{S})=\mathcal{S}\big] \ \mbox{iff} \ \
\big[\int_{(a,b]}|h(t)|^p |dr(t)|<\infty, \ \ \mbox{for all} \ \
0<p\le 2\big]
\]
\end{cor}
\begin{rem}
Taking on the unit interval $(0,1]$ the function $h(t)=t$ and the
time change $r(t):=t^{-\be}, \be
>2$, we see that the above corollary is not true for the mapping $I^{t,\,t^{-\be}}_{(0,1]}$
and all $0<p\le 2$.
\end{rem}
\medskip
\medskip
\textbf{7. Factorization property of measures from
$\mathcal{R}^{h,r}_{(a,b]}$}

\medskip
\emph{7.1. \underline{A motivating  example}}

Let us recall that for $\mathbb{B}_{t}=(B^{1}_t,B^{2}_t)$, Brownian
motion on $\mathbb{R}^{2}$, the process
\[
\mathcal{A}_{t} = \int_0^t B^{1}_sdB^{2}_s-B^{2}_sdB^{1}_s,\ \ t>0,
\]
is called \textit{L\'evy's stochastic area integral}. It is
well-known that for fixed $u>0$, and $a := (\sqrt{u},\sqrt{u})\in
\Rset^2 $ we have
\begin{equation*}
\chi(t) = E[e^{it\mathcal{A}_{u}}|\mathbb{B}_{u}= a] =
\frac{tu}{\sinh tu} \cdot \exp \{-(tu\coth tu-1)\},\quad t\in \Rset
\end{equation*}
(cf. L\'evy (1951) or Yor (1992), p. 19). If $\mu, \lambda$ and
$\nu$ are probability measures corresponding to the characteristic
functions $\chi(t), \phi(t): = tu/\sinh tu$ and $\psi(t): =
\exp[-(tu \coth tu-1)]$, respectively,  then $\lambda
=I^{e^{-t},t}_{(0,\infty)}(\nu)$, (cf. Remark 5(c)) and moreover
\[
\mu= I^{e^{-t},t}_{(0,\infty)}(\nu)\ast \nu \in \, L\equiv
I^{e^{-t},t}_{(0,\infty)}(ID_{\log})
\]
In other words, the conditional Levy's stochastic area integral has
selfdecomposable probability distribution $\mu$ that can be
\emph{factorized} by another selfdecomposable measure $\lambda$ and
its background driving measure $\nu$; cf. Iksanov, Jurek and
Schreiber (2004), p. 1367. That phenomena prompted the introduction
of the notion of factorization property for the L\'evy  class L
distributions.

\medskip
\emph{7.2. \underline{Definition and a condition for the
factorization property}}

\medskip
If for $\mu=I^{h,r}_{(a,b]}(\nu)\in\mathcal{R}^{h,r}_{(a,b]}$ we
also have that $I^{h,r}_{(a,b]}(\nu)\ast\nu\in
\mathcal{R}^{h,r}_{(a,b]}$ then we say that $\mu$ has \emph{a
factorization property}.
\begin{prop}
Suppose that for a given functions $h,r$ and an interval $(a,b]$
there exist function $h^{\prime},r^{\prime}$ and an interval
$(a^{\prime},b^{\prime}]$ such that for positive measures $\rho$ and
$\rho^{\prime}$, induced by the monotone functions $r$ and
$r^{\prime}$ respectively, we have the following
\begin{multline}
h((a,b])\cdot h^{\prime}((a^{\prime},b^{\prime}])=
h^{\prime}((a^{\prime},b^{\prime}])= h((a,b])=(c,d], \  \mbox{for some} \ 0<c<d, \\
\ \ \mbox{and} \ \ (h\otimes h^{\prime})(\rho \times \rho^{\prime})=
(h\rho)-(h^{\prime}\rho^{\prime})\ge 0. \qquad
\end{multline}
Then $\mathcal{D}^{h,r}_{(a,b]} \subset
\mathcal{D}^{h^{\prime},r^{\prime}}_{(a^{\prime},b^{\prime}]}$ and
for all $\nu\in\mathcal{D}^{h,r}_{(a,b]}$ putting
$\lambda:=I^{h^{\prime},r^{\prime}}_{(a^{\prime},b^{\prime}]}(\nu)$
we have
\[
I^{h^{\prime},r^{\prime}}_{(a^{\prime},b^{\prime}]}\big(I^{h,r}_{(a,b]}(\nu)\ast\nu\big)=I^{h,r}_{(a,b]}(\lambda)\ast
\lambda=I^{h,r}_{(a,b]}(\nu).
\]
In other words,
$\mathcal{R}^{h,r}_{(a,b]}=\{I^{h,r}_{(a,b]}(\lambda)\ast\lambda:
\lambda \in
I^{h^{\prime},r^{\prime}}_{(a^{\prime},b^{\prime}]}(\mathcal{D}^{h,r}_{(a,b]})
\}$
\end{prop}
\emph{Proof.} Since $0\le h^{\prime}\rho^{\prime} \le h\rho$
 then from Corollary 8 (expressed in terms of measures) we infer the                                                                                
inclusion of the domains.

From (23), Theorem 2 and the formula (3) we get
\begin{multline*}
\big(I^{h,\rho}_{(a,b]}\circ
I^{h^{\prime},\rho^{\prime}}_{(a^{\prime},b^{\prime}]}\big)(\nu)
\ast I^{h^{\prime},\rho^{\prime}}_{(a^{\prime},b^{\prime}]}(\nu) =
I^{h,\rho}_{(a,b]}(
I^{h^{\prime},\rho^{\prime}}_{(a^{\prime},b^{\prime}]}(\nu)) \ast
I^{h^{\prime},\rho^{\prime}}_{(a^{\prime},b^{\prime}]}(\nu)\\
=I^{t, (h\otimes h^{\prime})(\rho \times
\rho^{\prime})}_{(c,d]}(\nu) \ast I^{t,
(h^{\prime}\rho^{\prime})}_{(c,d]}(\nu)= I^{t,
(h\rho)}_{(c,d]}(\nu)= I^{h,\rho}_{(a,b]}(\nu),
\end{multline*}
which completes the proof.

The factorization property of a selfdecomposable measure given by
the Levy's stochastic area integral is not an exception as we have
\begin{cor}
For the class L of selfdecomposable probability measures on $E$ we
have
\[
L=\{I^{e^{-t}, t}_{(0,\infty)}(\nu)\ast\nu : \nu \in
I^{s,s}_{(0,1]}(ID_{\log}(E))\}
\]
\end{cor}
\emph{Proof.} We have that $L=I^{e^{-t},
t}_{(0,\infty)}(ID_{\log})$; cf. Jurek and Vervaat (1983). Then
taking $h^{\prime}(s)=s,\rho^{\prime}=\emph{l}_1$ (Lebesgue measure
on unit interval), $a^{\prime}=0$ and $b^{\prime}=1$ we check that
conditions (23) are fulfilled. Thus Proposition 8 gives the claim in
the corollary.

[The above fact was also shown in Jurek (2008), Theorem 3.1 but by a
different reasoning.]

\medskip
\medskip
\textbf{8. Some explicit examples.}

\medskip
\emph{8.1. \underline{Examples of domains of random integral
mappings}}

Here we recall a few examples of domains and in some instances
sketch their proofs that rely on Corollary 5.

\medskip
\begin{ex}
\begin{equation}
\mathcal{D}^{t, -\log t}_{(0,1]}=ID_{\log}(H):=\{\mu\in ID: \int_H
\log(1+||x||)M(dx)<\infty \}.
\end{equation}
\end{ex}
For this let us note that
\begin{multline*}
\int_H (1\wedge ||x||^2)\,M^{t, -\log
t}_{(0,1]}(dx)=\int_0^1\int_H(1\wedge
t^2||x||^2)\,M(dx)\,\frac{dt}{t}\\= \int_0^1 t \int_{(||x||\le
t^{-1})}||x||^2M(dx)dt +\int_0^1 \int_{(||x||> t^{-1})}
M(dx)\frac{dt}{t}= \\ 1/2\int_H (1\wedge ||x||^2)M(dx) +
\int_{(||x||>1)}\log ||x|| M(dx)<\infty,
\end{multline*}
which is equivalent with finite log-moment of $\mu$; cf. Jurek and
Smalara (1981) or Proposition 1.8.13 in Jurek and Mason (1993).

[Example 1 is valid on any Banach space E.  However, its proof is
completely different from the above for Hilbert space H; cf. Jurek
and Vervaat (1983).]

\begin{ex}

(1) $\mathcal{D}^{t, t^{\be}}_{(0,1]}=ID(E), \ \ \mbox{for} \ \
\be>0$.

(2) $\mathcal{D}^{t, t^{\be}}_{(0,1]}=ID_{\be}(H):= \{\nu\in
ID(H):\int_H||x||^{-\be}\nu(dx)<\infty \}$,

for  $-1<\be<0$.

(3)  $\mathcal{D}^{t,t^{\be}}_{(0,1]}\cap ID^{\circ}=ID_{\be}(H)\cap
ID^{\circ}, \ \  \mbox{for} \ \ \ -2< \be \le -1$;  where
$ID^{\circ}$ denotes symmetric infinitely divisible measures.
\end{ex}

\begin{rem}
Recall that the integral mappings $I^{t,t^{\be}}_{(0,1]}$ and their
domains appeared in the context of the classes $\mathcal{U}_{\be}$
for $-2\le \be <0$ and $0\le \be <\infty$. The class $\mathcal{U}_0$
coincides with the L\'evy class
$L=I^{e^{-t},t}_{(0,\infty)}(ID_{\log})$, while
$\mathcal{U}_{-2}(E)$ consists only od Gaussian measures; cf. Jurek
(1988), (1989)  and Jurek-Schreiber (1992).
\end{rem}

\medskip
We complete this subsection with examples of time changes given by
\emph{the incomplete Euler function.} It is defined as follows
\begin{equation}
\Gamma (\alpha; x):= \int_x^{\infty}t^{\alpha -1}\,e^{-t}dt, \ \ \
\alpha\in\Rset, \ \  x>0.
\end{equation}

For $\al>0$ the above is just the gamma function  and
$\Gamma(0+)<\infty$ and thus from Proposition 3 we get
\begin{ex}
For $ \al>0$, $\mathcal{D}^{t, \Gamma(\al;t)}_{(0,\infty)}=ID(E)$.
\end{ex}

Further for $\al=0$ we get
\begin{ex}
\begin{equation}
\mathcal{D}^{t,\,\,
\int_t^{\infty}\frac{e^{-s}}{s}{ds}}_{(0,\infty)}=ID_{\log}(H)
\end{equation}
\end{ex}
Similarly as in Example 1,
\begin{multline}
\int_0^{\infty}\int_H(1\wedge t^2||x||^2)\,M(dx)\,\frac{e^{-t}}{t}
dt \\= \int_0^{\infty} t \int_{(||x||\le t^{-1})}||x||^2M(dx)
e^{-t}dt
+\int_0^{\infty} \int_{(||x||> t^{-1})} M(dx)\frac{e^{-t}}{t}dt= \\
\int_{(||x||>0)} ||x||^2\, \big[\int_0^{||x||^{-1}} t e^{-t}dt
\big]\, M(dx) + \int_{(||x||>0)}\, \big[\int_{||x||^{-1}}^{\infty}
\frac{e^{-t}}{t} dt\big]\, M(dx).
\end{multline}
Note that in the first square bracket we get
\[
\int_0^{||x||^{-1}} t e^{-t}dt= 1- e^{-||x||^{-1}}(1+||x||^{-1})\le
1\wedge ||x||^{-2}
\]
and hence the first integral in (27) is finite.

For second integral in (27) let us brake the space
$H\setminus{\{0\}}$ into two parts. For $0<||x||\le 1$,
\begin{multline*}
\int_{(0<||x||\le 1)}||x||^2 \ [\
||x||^{-2}\int_{||x||^{-1}}^{\infty}
\frac{e^{-t}}{t} dt \ ] M(dx)\\
 \le [ \ \sup_{(a \ge
1)}\,a^2\int_a^{\infty}\frac{e^{-t}}{t} dt \ ]\int_{(||x||\le
1)}||x||^2 M(dx)< \infty.
\end{multline*}
 For the part $||x||>1$ we use Remark 4(a) that gives
\[
\int_{||x||^{-1}}^{\infty} \frac{e^{-t}}{t} dt = -\mathbf{C}+ \log
||x||+ \int_0^{||x||^{-1}}\frac{1-e^{-t}}{t}dt,
\]
where the integral on the right hand side is bounded by $\int_0^1
(1-e^{-t})t^{-1}dt<\infty$. All in all the second integral in (27)
is finite if and only if $\int_{(||x||>1)}\log ||x|| M(dx)<\infty$,
which completes the proof of Example 4.
\begin{ex}
For $-1<\al < 0$ we have
\[
\mathcal{D}^{t,\int
_t^{\infty}s^{\al-1}e^{-s}ds}_{(0,\infty)}=ID_{\al}(\Rset),
\]
\end{ex}
\noindent with the notations from Example 2. For the above example
and the case $-2<\al\le -1$ cf. Sato (2006).

\medskip
\medskip
\emph{8.2. \underline{Examples of iterated integral mappings and
image measures}}

\medskip
\begin{ex} For $\nu\in ID_{\log^{m}}$ and $m=1,2,...$
\[
I^{e^{-s},\,s}_{(0,\infty)}(I^{e^{-s},\,s}_{(0,\infty)}(...
I^{e^{-s},\,s}_{(0,\infty)}(\nu)))=
I^{e^{-t},\,\frac{t^{m}}{m!}}_{(0,\infty)} (\nu)
\]
\end{ex}
\emph{Proof.} In view of Theorem 2 it is enough to check that for
$h(t):=e^{-t}$ and $\rho:=\emph{l}$ (the Lebesque measure on
$\Rset$) we have equality
\begin{multline*}
\int_0^{\infty}g(u)[(e^{-t})^{\otimes m}) (\emph{l}_1\times ...
\times \emph{l}_1)](du) \\=\int_0^{\infty} \int_0^{\infty}...
 \int_0^{\infty}g(e^{-(s_1+s_2+...+s_m)}u)ds_1...ds_m= \int_0^{\infty}g(e^{-s}u)d\,[\frac{s^m}{m!}]
\end{multline*}
for all $g$ bounded and measurable. The first equality is just a
change of variable argument. For the second, using the induction
arguments, we have
\begin{multline*}
\int_0^{\infty} \Big[\int_0^{\infty}...
 \int_0^{\infty}g(e^{-(s_1+s_2+...+s_{m-1})}e^{-s_m}u)ds_1...ds_{m-1}\Big]ds_m =\\
 \int_0^{\infty}\int_0^{\infty}g(e^{-t}
 e^{-s_m}u)d\,[\frac{t^{m-1}}{(m-1)!}]ds_m
 =\int_0^{\infty}\int_{s_m}^{\infty}g(e^{-w}u)\frac{(w-s_m)^{m-2}}{(m-2)!}dwds_m
 \\ =
 \int_0^{\infty}g(e^{-w}u)\frac{w^{m-1}}{(m-1)!}dw =
 \int_0^{\infty}g(e^{-w}u)d\,[\frac{w^{m}}{m!}],
\end{multline*}
which completes the proof.

\begin{rem} The class of measures
$I^{e^{-t},\,\frac{t^{m}}{m!}}_{(0,\infty)} (ID_{\log^m})$ coincides
with the set $L_m$ of so called \emph{m-times selfdecomposable
distributions}; cf. Jurek (2011) for the history of those classes
and relevant  references.
\end{rem}
\begin{ex} For $\beta>0$ we have
\begin{equation}
I^{t^{1/\beta},\, t}_{(0,1]}\circ I^{s^{1/2\beta},\,
s}_{(0,1]}=I^{w,\, 2w^{\beta}(1-(1/2)
w^{\beta})}_{(0,1]}=I^{(1-\sqrt{t})^{1/\beta},\, t}_{(0,1]}
\end{equation}
Or equivalently, for Lebesque measure $l_1$ on the unit interval and
$ 0< w \le 1$ we get
\begin{equation*}
(t^{1/\beta}\otimes s^{1/(2\beta)})(l_1 \times l_1)(dw)= id
^{\otimes 2}(\beta t^{\beta-1}dt \times 2 \beta s^{2
\beta-1}dt)(dw)=2 \beta w^{\beta-1}(1-w^{\beta})\,dw
\end{equation*}
\end{ex}
\emph{Proof.} As in Example 6, it simply follows from Theorem 2 and
identity (3) because all time change functions are strictly
increasing on the unit interval.
\begin{ex}
For $\beta >0$
\[
I^{t^{1/\beta},\, t}_{(0,1]}\circ I^{e^{-s},\,s}_{(0,\infty)}=
I^{e^{-s},\,\,s+\beta^{-1}e^{-\beta
s}-\beta^{-1}}_{(0,\infty)}=I^{-w,\,\, \beta^{-1}w^{\beta}-\log
w-\be^{-1}}_{(0,1]}.
\]
Or equivalently, for $0<w\le1$
\begin{equation*}
(t^{1/\beta}\otimes e^{-s})(l_1 \times l)(dw)=
(\beta^{-1}w^{\beta}-\log w-\be^{-1})dw.
\end{equation*}
\end{ex}
\noindent This is a consequence of Theorem 2. Also cf. Czyżewska-Jankowska
and Jurek (2011), Proposition 2.
\begin{ex} For $\al\in \Rset$
\[
I^{t,\, \Gamma(\al; t)}_{(0,\infty)}\circ
I^{e^{-s},\,s}_{(0,\infty)}=
I^{t,\,\int_t^{\infty}s^{-1}\Gamma(\al;s)ds}_{(0,\infty)},
\]
\end{ex}
\noindent which follows from Theorem 2.

\medskip
\noindent Institute of Mathematics, University of Wroc\l aw, Pl.
Grunwaldzki 2/4, 50-384 Wroc\l aw, Poland. [E-mail:
zjjurek@math.uni.wroc.pl]


\begin{thebibliography}{}
\bibitem{} T. Aoyama (2009), Nested subclasses of the class of type G
selfdecomposable distributions on $\Rset^d$, \emph{Probab. Math.
Stat.}, vol. 29,  pp. 135-154.

\bibitem{AG}  A. Araujo, and E. Gin\'e (1980), \textit{The Central Limit Theorem for Real
and Banach Valued Random Variables}, Wiley, New York.


\bibitem{} O. E. Barndorff-Nielsen, M. Maejima and K. Sato (2006),
Some classes of multivariate infinitely divisible distributions
admitting integral representations, \emph{Bernoulli} \textbf{12},
pp. 1-33.

\bibitem{} P. Billingsley (1968), \emph{Convergence of probability
measures}, Wiley, New York.


\bibitem{} A. Czyżewska-Jankowska and Z. J. Jurek (2011),
Factorization property of generalized s-selfdecomposable measures
and class $L^f$ distributions, \emph{Theory Probab. Appl.} vol.
\textbf{55}, no 4, pp. 692-698.

\bibitem{}  I. S. Gradshteyn and I. M. Ryzhik (1965).
\emph{Tables of integrals, series, and products}, Academic Press,
New York.

\bibitem{}  B. Grigelionis (2007), Extended Thorin classes and
stochastic integrals, \emph{Liet. Matem. Rink.} \textbf{47} , pp.
497 -- 503.

\bibitem{} A.M. Iksanov, Z. J. Jurek and B.M. Schreiber (2004), A new
factorization property of the selfdecomposable probability measures,
\textit{Ann. Probab.} \textbf{32}, no 2, 1356--1369.


\bibitem{} Z. J. Jurek (1983), Limit distributions and one-parameter groups of linear operators on Banach spaces,
\emph{J. Multivar. Anal.} \textbf{13}, pp. 578-604.

\bibitem{J88} Z.J. Jurek (1988), Random integral representations for
classes of limit distributions similar to Levy class $L_0$,
\textit{Probab. Th. Fields} \textbf{78}, 473-490.

\bibitem{} Z. J. Jurek (1989), Random integral representations for classes of limit
distributions similar to L\'evy class $L_{0}$. II. \emph{Nagoya
Math. Journal} \textbf{114}, pp. 53-64.

\bibitem{07} Z.J. Jurek (2007), Random integral representations for free-infinitely divisible
and tempered stable distributions, \emph{Stat.\& Probab. Letters},
\textbf{77} no. 4, pp. 417-425.


\bibitem{} Z. J. Jurek (2008), A calculus on
L\'evy exponents and selfdecomposability on Banach spaces,
\emph{Probab. Math. Stat.} vol. \textbf{28}, Fasc. 2, pp. 271-280.

\bibitem{} Z. J. Jurek (2011), The Random Integral Representation Conjecture:  a quarter of a century later,
\emph{Lithuanian Math. Journal}, \textbf{51}, no 3, 2011, pp.
362-369.

\bibitem{JM93} Z. J. Jurek and J.D. Mason (1993), \textit{Operator-limit Distributions
in Probability Theory}, Wiley Series in Probability and Mathematical
Statistics, New York.

\bibitem{} Z. J. Jurek and B.M. Schreiber (1992), Fourier transforms of measures from classes $U_{\beta},-2 <\beta
<-1$.  \emph{J. Multivariate Analysis} \textbf{41} (1992), pp.
194-211.

\bibitem{} Z. J. Jurek and J. Smalara (1981), On integrability with
respect to infinitely divisible measures, \emph{Bull. Acad. Polon.
Sci.} \textbf{29}, pp. 179-185.

\bibitem{JV83} Z.J. Jurek and W. Vervaat (1983), An integral representation for
selfdecomposable Banach space valued random variables, \textit{Z.
Wahrsch. verw. Gebiete} \textbf{62}, 247--262.

\bibitem{L} P. L\'evy (1951) , Wiener's random functions, and other
Laplacian random functions. In \textit{Proc. 2nd Berkeley Symposium
Math. Stat. Probab.}, Univ. California Press, Berkeley, 171--178.

\bibitem{} W. Linde (1983), \emph{Probability measures in Banach spaces- stable and
infinitely divisible distributions}, Wiley, New York.

\bibitem{} M. Maejima, V.Perez Abreu, K.I. Sato (2012), A class of
multivariate infinitely divisble distributions related to arcsine
density, \emph{Bernoulli}, vol. 18 no. 2, pp. 476-495.

\bibitem{} M. Maejima and K. Sato (2009), The limits of nested subclasses of
several classes of infinitely divisible distributions are identical
with closure of the class od stable distributions, \emph{Probab.
Rel. Fields}, vol. 145, pp. 119-142.

\bibitem{P} K. R. Parthasarathy (1968), \textit{Probability measures on
metric spaces}, Academic Press, New York and London, 1968.

\bibitem{} K. Sato (1999),\emph{L\'evy processes and infiniteley divisible
distributions}, University Press, Cambridge, United Kingdom.

\bibitem{} K. Sato (2006), Two families of improper stochastic
integrals with respect to L\'evy processes, \emph{ALEA, Lat. Am. J.
Probab. Math. Stat.} vol. 1, pp. 47-87.

\bibitem{Y} M. Yor (1992), \textit{Some Aspects of Brownian Motion, Part I:
Some Special Functionals}, Birkh\"auser, Basel.


\end{thebibliography}
\end{document}